\documentstyle{article}

                    \newcommand{\w}{\widetilde}
                    \renewcommand{\l}{\lambda}
                    \newcommand{\fr}{\frac}

                    \newcommand{\bequ}{\begin{equation}}
                    \newcommand{\eequ}{\end{equation}}
                    \newcommand{\beqa}{\begin{eqnarray}}
                    \newcommand{\eeqa}{\end{eqnarray}}
                    \newcommand{\beq}{\begin{eqnarray*}}
                    \newcommand{\eeq}{\end{eqnarray*}}
										\newcommand{\barr}{\begin{array}}
\newcommand{\earr}{\end{array}}
\renewcommand{\O}{{\cal O}}
\newcommand{\bfr}{\begin{flushright}}
\newcommand{\efr}{\end{flushright}}
\newcommand{\bfl}{\begin{flushleft}}
\newcommand{\efl}{\end{flushleft}}
\renewcommand{\O}{{\cal O}}
										
\newtheorem{opr}{Definition}
\newtheorem{sle}{Corollary}
\newtheorem{teo}{Theorem}
\newtheorem{lem}{Lemma}
\newtheorem{pre}{Supposition}
\newcommand{\bo}{\begin{opr}}
\newcommand{\eo}{\end{opr}}
\newcommand{\bs}{\begin{sle}}
\newcommand{\es}{\end{sle}}
\newcommand{\bt}{\begin{teo}}
\newcommand{\et}{\end{teo}}
\newcommand{\bl}{\begin{lem}}
\newcommand{\el}{\end{lem}}
\newcommand{\bp}{\begin{pre}}
\newcommand{\ep}{\end{pre}}
\newcounter{rem}
\newcounter{pr}
\newcounter{remar}
\newcounter{cpr}

\font\Sets=msbm10

\def\Re{\hbox{\Sets R}}

\title{An analogue of Borg's uniqueness theorem
in the case of indecomposable boundary conditions}
\author{Azamat M. Akhtyamov}

\begin{document}
\maketitle

\begin{abstract}
An uniqueness theorem for
 the inverse problem in the case of a second-order
equation defined on the interval [0,1]
when the boundary forms contain
combinations of the values of functions
at the points 0 and 1 is proved.
The auxiliary
eigenvalue problems in our theorem
 are chose
in the same manner as in Borg's  uniqueness theorem are not
as in that of Sadovni\v ci$\check \imath $'s.
So number of conditions in our theorem is
less than that in  Sadovni\v ci$\check\imath$'s.
\end{abstract}

In the case of indecomposable boundary conditions
(i.~e. Sturm-Liouville conditions)
several methods for the solution of the inverse problem
are worked out
in \cite{Ambarzumijan 29}--\cite{Hochstadt 73}.
Extensive bibliographies for this problem can
be found in  Levitan  \cite{Levitan 84}.
In \cite{Sadovnichiy 72}
a new method
 was proposed for the solution of
the inverse problem with indecomposable boundary conditions.
In this method instead of trasformation
operators are utilized mappings $T_{\l }$ of
the spaces of the solutions defined by
matrices.

It was shown that to restore uniquely the
function~$q(x)$ and the boundary conditions
by a
 set of the eigenvalues $\l _k$ of
eigenvalue problem itself,
a sets of the eigenvalues $z_{k,\, a}$ of two auxiliary
eigenvalue problems,
the ''weight'' numbers $\alpha _{k,\, a}$ and resides of the
certain functions.

Sadovni\v ci$\check\imath$  showed
the ''weight'' numbers $\alpha _{k,\, a}$ and resides of the
certain functions
is necessary
for the unique restoration
of the function $q(x)$ and the boundary conditions.

 We will show that if
two auxiliary
eigenvalue problems
 are chose
in the same manner as in Borg's  uniqueness theorem (\cite{Borg 46}) are
not
as in that of Sadovni\v ci$\check \imath $'s (\cite{Sadovnichiy 72}),
then
 the ''weight'' numbers $\alpha _{k,\, a}$ and resides of the
certain functions
are
unnecessary for the
 uniqueness theorem.

\vspace{0.2cm}

{\bf 1.} We will consider the following problem:
\beqa
ly &=& - y'' + q(x)\, y = \l \, y ,
\label{lanl-s 1}
\\
U_1(y) &=& y'(0) + a_{11}\, y(0) + a_{12}\, y(1) = 0,
\label{lanl-s 2}
\\
U_2(y) &=& y'(1) + a_{21}\, y(0) + a_{22}\, y(1) = 0
\label{lanl-s 3}
\eeqa
($a_{ij},$ $i,\, j =1,\, 2$ are
real constants, $q(x) \in C^1[0,\, 1].$).

The problem defined in (\ref{lanl-s 1})--(\ref{lanl-s 3})
will be called problem $L.$
We denote by $\w L$ a problem of type $L$ but with different
coefficients in the equation and with different parameters
 in the boundary forms.
In all that follows, if certain symbol denotes a term from problem $L,$
then the symbol $\w {}$ denotes the analogous term from problem
$\w L.$ Everywhere the integral index $k$ varies from $0$ to $\infty .$

Along with problem $L$ we consider two problems with
decomposable
boundary conditions:
problems

$L_1:$
\beq
ly &=& - y'' + q(x)\, y = \l \, y ,
\\
U_{1,1}(y) &=& y'(0) + a_{11}\, y(0) = 0,
\\
U_{2,1}(y) &=& y'(1) + a_{22}\, y(1) = 0
\eeq

and $L_2:$
\beq
ly &=& - y'' + q(x)\, y = \l \, y ,
\\
U_{1,2}(y) &=& y'(0) + a_{11}\, y(0) = 0,
\\
U_{2,2}(y) &=& y'(1) + a_{12}\, y(1) = 0
\eeq

Let $\l _{k,\, 1}$, $\l _{k,\, 2}$
be eigenvalues of these problems.

\bt
Let
$\l _{k}=\w \l _{k} ,$
$\l _{k,\, 1}=\w \l _{k,\, 1} ,$
$\l _{k,\, 2}=\w \l _{k,\, 2} ,$
$a_{12}\ne a_{22};$
then the coefficients of the equations and the constants
in the boundary conditions of the problems
$L$ and $\w L$ coincide, i.e.
$q(x)=\w q(x),$ $a_{ij}=\w a_{ij}$ $i,\, j =1,\, 2.$
\et

{\sl Proof.}  Making use of
Borg's  uniqueness theorem (\cite{Borg 46,Levitan 84}) to
problems $L_1$ and $L_2$ we have:
\bequ
q(x)=\w q(x),\qquad a_{11}=\w a_{11}, \qquad a_{12}=\w a_{12},
\qquad
a_{22}=\w a_{22}.
\label{lanl-s =1}
\eequ

To prove the theorem we must show that $a_{21}=\w a_{21}.$
We shall do it.

Let $y_1(x,\, \l )$ and $y_2(x,\, \l )$ be
linearly independent solutions of
equation (\ref{lanl-s 1}) satisfying
$$y_1(0,\, \l )=1,\quad y_1'(0,\, \l )=0,\quad y_2(0,\, \l )=0,\quad
y_2'(0,\, \l )=1.$$
Then
asymptotic formulae
\beq
y_1(x,\, \l) &=& \cos \l x + \fr{1}{\l }\, u(x)\, \sin \l x +\O \left(
\fr{1}{\l
^2}\right) ,
\\
y_2(x,\, \l) &=& \fr{1}{\l }\, \sin \l x - \fr{1}{\l ^2}\, u(x)\, \cos \l x
+\O \left(
\fr{1}{\l
^3}\right) ,
\\
y_1'(x,\, \l) &=& {\l }\, \sin \l x + u(x)\, \cos \l x +\O \left(
\fr{1}{\l
}\right) ,
\\
y_2'(x,\, \l) &=& \cos \l x + \fr{1}{\l }\, u(x)\, \sin \l x +\O \left(
\fr{1}{\l
^2}\right)
\\
\mbox{where } u(x) &=&  \fr{1}{2}\, \int_0^x q(t) \, dt
\eeq
are true
for $\l \in \Re $ and  $\l $ sufficiently large
 (\cite{Naymark 68}).

 It follows according to the condition of the
theorem, that the eigenvalues
of spectral problems
$L$ and $\w L$
coincide.
The eigenvalues of the problem
$L$
are the roots of the following entire function of the first order
$$
\Delta(\l )
=
\left|
\barr{cc}
U_1(y_1(x,\l )) &
U_1(y_2(x,\l )) \\
U_2(y_1(x,\l )) &
U_2(y_2(x,\l ))
\earr
\right|
$$
(\cite{Naymark 68}).

If we substitute
the asymptotic formulae
of solutions
$y_1(x,\l ),$
$y_2(x,\l ))$
 for
the preceding equation, we obtain
\beq
\Delta(\l )
&=&
(a_{11}+a_{12}\, y_1(1,\, \l ))\cdot (y_2'(1,\, \l )
+ a_{22}\, y_2(1,\, \l )) -
\\
&&- (1 + a_{12}\, y_2(1,\, \l ))\cdot (y_1'(1,\, \l )
+a_{21} + a_{22}\, y_1(1,\, \l )) =
\\
&&\hspace{-1.5cm}=
a_{11}\,\cos \l + a_{12} + \l\,\sin \l - u_1(1)\, \cos ^2\l
+ a_{21} + a_{22}\, \cos \l + \O \left( \fr{1}{\l }\right).
\eeq

It follows from Weierstrass theorem
about an entire function representation by its roots that
$$
\Delta (\l )
\equiv
e^{a\l + b}
\w \Delta(\l ),
$$
where $\widetilde\Delta (\l )$ is
a characteristic determinant
of the problem $\w L$ and
$a,$ $b$ are  certain numbers.

From here
\beq
&&\widetilde\Delta (\l )
-
e^{a\l + b}
\Delta(\l )
\; \equiv \;
\\
&&\hspace{0.5cm}\; \equiv \;
(a_{11} - \w a_{11}\, e^{a\l + b})\,\cos \l
+ (a_{12} - \w a_{12}\, e^{a\l + b} ) +
\\
&&\hspace{1cm}+ (1 - e^{a\l + b})\,\l\,\sin \l
-(1-e^{a\l + b})\, u_1(1)\, \cos \l  -
\\
&&\hspace{1cm}- (a_{21} - \w a_{21}\, e^{a\l + b}) -
(a_{22}- \w a_{22}\, e^{a\l + b})\, \cos \l +
\\
&&\hspace{1cm}+
(1- e^{a\l + b})\,\O \left( \fr{1}{\l }\right)
\;\equiv \; 0.
\eeq

The functions
$1,$  $\sin \l ,$   $\cos \l ,$
 $\cos 2\l ,$ $\l\cdot \sin \l ,$
$\O \left( \fr{1}{\l }\right)$
 are linear independent functions with
respect to argument $\l $.
(It is easily verified by the definition of functions linear
independence.)
Consequently
$a=0,$ $b=0$ and
$$
(a_{12} - \w a_{12}) +
(a_{21} - \w a_{21}) +
(a_{22}- \w a_{22})\, \cos \l +
\O \left( \fr{1}{\l }\right) \equiv
 0.
$$

This together with (\ref{lanl-s =1}) then gives
$a_{21} = \w a_{21}.$ \ $\Box $


\begin{thebibliography}{99}


\bibitem{Ambarzumijan 29}
V.~A.~Ambarzumijan,
 {Uber eine Frage der Eigenwerttheorie} //
 Zeitshrift fur Physik,
 53 (1929), 690--695

\bibitem{Borg 46}
G.~Borg, {Eine Umkehrunrung der Sturm--Liouvillschen
  Eigenwertanfgabe. Bestimmung der Differentialgleichung
  durch die Eigenwarte} //
 Acta Mathematica, 78 (1946), N~2, 1--96.

\bibitem{Levinson 49}  N.~Levinson, The inverse Sturm--Liouville problem //
            Math. Tidsskr., Ser.B., 13 (1949), 25--30.

\bibitem{Gelfand 51} I.~M.~Gelfand and B.~M.~Levitan,
On the determination of a differential equation
from its spectral function,
Izv. Akad. Nauk SSSR Ser. Mat. 15 (1951), 309--360;
English transl., Amer. Math Soc. Transl. (2) 1 (1955), 253--304.

\bibitem{Krein 51}
M.~G.~Krein, Solution of the inverse Sturm-Liouville
problem, Dokl. Akad. Nauk SSSR 76 (1951), 21--24. [in Russian]
MR 12, 613.

\bibitem{Marchenko  52}
V.~A.~Marchenko, Some questions of the theory of
one-dimensional linear differential operators of the second order,
Trudy Moskov. Mat. Ob\v s\v c.

\bibitem{Levitan 84} B.~M.~Levitan,
            The inverse Sturm--Liouville problems and
                  applications [in Russian]
           Nauka, Moscow.
            1984.

\bibitem{Leybazon 66} Z.~L.~Leybazon,
            An inverse problem of the spectral analysis
                  of order differential operators of
higher order [in Russian]
            Trudy. Moskov. Mat. Ob\v s\v c. 1966. V.~15. P.~70--144 =
Trans. Moscow Math. Soc. 1966, 78--163. MR~34 \# 4951.

\bibitem{Hochstadt 73} H.~Hochstadt,
            The inverse Sturm--Lioville problem //
            Comm. Pure Appl. Math.
					   26 (1973), 715--729


\bibitem{Sadovnichiy 72} V.A.~Sadovni\v ci\v i
Uniqueness of the solution to the inverse problem
in the case of a second-order equation with the indecomposable
                   boundary conditions, regularized sums
                   of some of the eigenvalues. Factorization of
                   the characteristic determinant. [in Russian]
           Dokl. Akad. Nauk SSSR. 1972.
           Tom 206. No.2. P.~293--296. =
Soviet Math. Dokl. Vol.~13 (1972), No.~5, 1220--1223.

\bibitem{Naymark 68}
  M.~A.~Naymark  Linear differential operators. II. Ungar. New York. 1968.


\end{thebibliography}
\end{document}